\documentclass[review]{elsarticle}
\usepackage{amsfonts,amsmath,upref,amssymb,amsthm,amsxtra,amscd,enumerate}
\usepackage{color}
\DeclareMathOperator{\diver}{div}

\newcommand{\dx}{\,\mathrm{d}x}
\newcommand{\dy}{\,\mathrm{d}y}

\newcommand{\dt}{\,\mathrm{d}t}

\def\XXint#1#2#3{{\setbox0=\hbox{$#1{#2#3}{\int}$ }
\vvcenter{\hbox{$#2#3$ }}\kern-.6\wd0}}

\newtheorem{theorem}{Theorem}[section]

\newtheorem{definition}[theorem]{Definition}

\newtheorem{remark}[theorem]{Remark}

\numberwithin{equation}{section}

\newcommand{\esssup}{\mbox{ess sup }}

\newcommand{\weak}{\mbox{weak}}

\newcommand{\vep}{\varepsilon}
\newcommand{\pat}{\partial_t}

\def\XXint#1#2#3{{\setbox0=\hbox{$#1{#2#3}{\int}$ }
\vvcenter{\hbox{$#2#3$ }}\kern-.6\wd0}}

\newtheoremstyle{dotless}{}{}{\itshape}{}{\bfseries}{}{ }{}

\theoremstyle{dotless}

\renewcommand{\u}{{\bf u}}
\newcommand{\vv}{{\bf v}}

\newcommand{\U}{{\bf U}}
\newcommand{\bphi}{{\boldsymbol \varphi}}

\begin{document}

\begin{frontmatter}

\title{Inviscid limit for the compressible Euler system with non-local interactions}
\author{Jan B\v rezina\fnref{honza}}
\address{Tokyo Institute of Technology,\\ 2-12-1 Ookayama, Meguro-ku,\\ Tokyo, 152-8550\\ Japan}
\fntext[honza]{email:brezina@math.titech.ac.jp}

\author{V\' aclav M\' acha\fnref{venca}}
\address{Industry-University Research Center, Yonsei University,\\ 50 Yonsei-ro Seodaemun-gu,\\ Seoul, 03722\\ Republic of Korea}
\fntext[venca]{email:macha@math.cas.cz\\ The research of V. M. has been supported by the grant NRF-20151009350. }

\begin{abstract}
The collective behavior of  animals can be modeled by a system of equations of continuum mechanics endowed with extra terms describing repulsive and attractive forces between the individuals. This system  can be viewed  as a generalization of the compressible Euler equations with all of its unpleasant consequences, e.g., the non-uniqueness of solutions. In this paper, we analyze the equations describing a viscous approximation of a generalized compressible Euler system and we show that its dissipative measure-valued solutions tend to a strong solution of the Euler system as viscosity tends to 0, provided that the strong solution exists.
\end{abstract}


\end{frontmatter}

\section{Introduction}
Recently the modeling of collective behavior of animals has became a popular topic. In this paper we study the hydrodynamical model for collective behavior which can be obtained as a mean-field limit of the Cucker-Smale flocking model (see \cite{CaCaRo, KaMeTr}). As the hydrodynamical model is in fact an Euler-type system, its solutions may be constructed by means of convex integration (see \cite{CaFeGwSw, ChDeKr, ChKr, FeKr}) and thus are in general non-unique and irregular even for smooth initial data. Then the question of choosing the unique physically reasonable  solution comes into play (see  \cite{CaFeGwSw}). The inviscid limit is often considered to be one of the criteria of admissibility for the Euler-type problems. In this paper we study the question whether a solution to the above mentioned Euler type system could be understood as a limit of viscous approximations. 

More precisely, we prove that under certain smoothness conditions imposed on a solution to the Euler system \eqref{zaklad}, this solution can be  understood as a limit of solutions to the approximate Navier-Stokes system \eqref{art.vis}. 

We consider the following system:
\begin{equation}
\begin{split}\label{zaklad}
\partial_t \rho + \diver_x(\rho \u) =\ & 0,\\
\partial_t(\rho \u) + \diver_x(\rho \u\otimes \u) + \nabla_x p(\rho) =\ & (1-H(|\u|^2))\rho\u\\ &- \rho \int_{\Omega} \nabla_x K(x-y)\rho(y) \dy \\&+ \rho \int_{\Omega} \psi(x-y) \rho(y) \left(\u(y)-\u(x)\right) \dy,
\end{split}
\end{equation}
where $\rho= \rho(t,x)$ and $\u = \u(t,x)$  denote the unknown  density and velocity; $p=p(\rho)$ is the pressure; $H = H(z)$, $z\in [0,+\infty)$ is 'friction'; $ K = K(x)$ and $\psi = \psi(x)$ represent the  non-local interaction forces acting on the medium.

 For the sake of simplicity, we consider the spatially periodic boundary conditions, i.e.,
\begin{equation*}
\Omega= ([-1,1]|_{-1,1})^N,\ \ N=2,3.
\end{equation*}
In accordance with \cite{CaFeGwSw} we assume  that
\begin{equation}\label{tlak.pred}
\begin{split} p\in C([0,\infty))\cap C^2(0,\infty), \quad p(0) &= 0,\quad p'(\rho) >0 \mbox{ whenever }\rho >0,\\
 \liminf_{\rho \to \infty} p'(\rho) >0,& \quad \liminf_{\rho \to \infty} \frac{P(\rho)}{p(\rho)}
 > 0,
\end{split}\end{equation}
where
$$
 P(\rho) = \rho\int_1^\rho \frac{p(z)}{z^2} {\rm d}z;
$$
together with
\begin{equation}\label{h.pred}H\in C^2([0,\infty)),\quad 0 \leq H(z) \leq H_\infty:= \lim_{z\to \infty} H(z), \quad H'(z) \geq 0 \mbox{ for all } z\geq z_0;\end{equation}
and
\begin{equation}\label{k.psi.pred}
\begin{split} 
K\in C^2(\Omega), \quad \psi&\in C^1(\Omega), \quad \psi\geq 0,\\ K \mbox{ and } \psi \mbox{  symmetric, i.e., } &K(x) = K(-x), \ \psi(x) = \psi(-x).
\end{split}
\end{equation}
Note that we could take $H,\ K$ and $\psi$ so that no extra terms appear on the right-hand side of \eqref{zaklad}$_2$. Therefore, our result also covers the case of  compressible Euler equations.
Finally, the system is endowed with the following initial conditions:
\begin{equation}\label{init.con}
\rho(0) = \rho_0,\ \quad \u(0) = \u_0.
\end{equation}

We show that any sufficiently smooth solution to \eqref{zaklad} can be obtained as a limit of dissipative measure-valued solutions to the viscous approximation 
\begin{equation}
\begin{split}\label{art.vis}
\partial_t \rho^\vep + \diver_x(\rho^\vep \u^\vep) =\ & 0,\\
\partial_t(\rho^\vep \u^\vep) + \diver_x(\rho^\vep \u^\vep\otimes \u^\vep) + \nabla_x p(\rho^\vep) =\ &\vep \Delta_x \u^\vep + (1-H(|\u^\vep|^2))\rho^\vep \u^\vep \\ - \rho^\vep& \int_{\Omega} \nabla_x K(x-y)\rho^\vep(y) \dy \\ + \rho^\vep& \int_{\Omega} \psi(x-y) \rho^\vep(y) \left(\u^\vep(y)-\u^\vep(x)\right) \dy,
\end{split}
\end{equation}
 on $\Omega$ as $\vep$ tends to zero. The unknowns $\rho_\vep$ and $\u_\vep$ are assumed to satisfy the initial conditions
\begin{equation}\label{art.vis.init.con}
\rho^\vep(0) = \rho^\vep_0, \ \quad \u^\vep(0) = \u^\vep_0,
\end{equation}
such that $\rho^\vep_0 \to \rho_0$ and $\u^\vep_0 \to \u_0$ in a sense specified later on.

As there are two crucial non-linearities in \eqref{art.vis}$_2$, namely $p(\rho^\vep)$ and $(1-H(|\u^\vep|^2))\rho^\vep \u^\vep$, it seems that the existence of weak solutions to  \eqref{art.vis} cannot be  established by the currently available theory. Thus we turn our attention towards the measure-valued solutions. Motivated by \cite{FeGwSwWi} we define a measure-valued solution to \eqref{art.vis} and \eqref{art.vis.init.con} in the following way.
\begin{definition} \label{dis.mes.sol}
Let $\varepsilon > 0$,  $T > 0$  and $Q_\tau:=(0,\tau)\times \Omega$ for every $\tau \in [0,T]$.  We say that a parametrized measure $\{\nu_{t,x}\}_{(t,x)\in Q_T}$,
$$
\nu\in L^\infty_{w*}(Q_T, \mathcal P([0,\infty)\times \mathbb R^N), \ \rho(t,x) := \langle \nu_{t,x}, s\rangle,\ \u(t,x) := \langle \nu_{t,x}, \vv \rangle,\footnote{Hereinafter we use convention $\langle \nu_{t,x},F(s,\vv)\rangle = \int_{[0,\infty)\times \mathbb R^N} F(s,\vv)\ \rm{d}\nu_{t,x}$ where $s$ and $\vv$ represent the dummy variables for $\rho$ and $\u$, respectively.}
$$
is a dissipative measure-valued solution to \eqref{art.vis} in $Q_T$  and \eqref{art.vis.init.con} in $\Omega$ with the dissipation defect $\mathcal D\in L^\infty(0,T)$, $\mathcal D\geq 0$ if:
\begin{itemize}
\item For every $\psi \in C^1([0,T]\times \Omega)$, the {\bf equation of continuity}
\begin{multline*}
\int_\Omega \rho(\tau)  \varphi(\tau)\dx - \int_{Q_\tau} \rho \pat \varphi \dx \dt \\- \int_{Q_\tau}\langle \nu_{t,x}, s\vv\rangle \nabla_x \varphi(t,x)\dx \dt = \int_\Omega \langle {\nu_{0}}_x,s \rangle\varphi (0,x) \dx,
\end{multline*}
holds for all  $\tau \in [0,T]$.

\item There exists a measure $ \mu\in L^\infty(0,T, \mathcal M(\Omega))$ and a constant $c > 0$ such that for every $\bphi \in C^1([0,T]\times \Omega)^N$, 
$$
\left|\int_\Omega \bphi(\tau) {\rm d}\mu(\tau)\right|\leq c \mathcal D(\tau) \|\bphi(\tau)\|_{L^\infty(\Omega)},
$$
and the {\bf momentum equation}
\begin{multline*}
\int_\Omega \langle \nu_{\tau,x}, s\vv\rangle \bphi(\tau,x) \dx  - \int_{Q_\tau}\langle \nu_{t,x},s\vv\rangle\pat \bphi(t,x)\dx \dt\\ - \int_{Q_\tau} \langle \nu_{t,x}, s\vv\otimes \vv\rangle \nabla_x \bphi(t,x)\dx \dt  + \varepsilon\int_{Q_\tau}\nabla_x \u:\nabla_x \bphi \dx \dt \\ - \int_{Q_\tau} \langle \nu_{t,x},p(s)\rangle \diver_x \bphi(t,x) \dx \dt \\ =\int_{Q_\tau} \langle \nu_{t,x},  (1-H(|\vv|^2))s\vv\rangle \bphi(t,x)\dx \dt \\- \int_{Q_\tau} \rho(t,x) \int_\Omega \nabla_xK(x-y) \rho(t,y) \dy \bphi(t,x)\dx \dt\\
 + \int_{Q_\tau} \rho(t,x) \int_\Omega \psi(x-y)\langle \nu_{t,y},s\vv \rangle\dy \bphi(t,x)\dx \dt
 \\- \int_{Q_\tau} \langle \nu_{t,x}, s\vv\rangle \int_\Omega \psi(x-y)\rho(t,y)\dy \bphi(t,x)\dx \dt\\
 + \int_\Omega \langle {\nu_0}_x, s\vv\rangle \bphi(0,x) \dx + \int_0^\tau \int_\Omega \nabla_x \bphi(t) {\rm d}\mu(t)\dt,
\end{multline*}
hold for all $\tau \in [0,T]$. 

\item The {\bf energy inequality} 
\begin{multline*}
\int_\Omega \left\langle\nu_{\tau,x}, \frac 12s|\vv|^2 + P(s) \right\rangle\dx + \frac 12\int_\Omega (\rho K*\rho)(\tau)\dx\\ + \vep \int_{Q_\tau}  |\nabla_x\u|^2 \dx\dt + \mathcal D(\tau) \\ 
\leq \int_{\Omega} \left\langle {\nu_0}_x,  \frac 12 s |\vv|^2 + P(s) + \frac 12 sK*\langle {\nu_0}_y, s \rangle\right\rangle \dx\\
+\int_{Q_\tau} \langle\nu_{t,x}, (1- H(|\vv|^2)) s|\vv|^2\rangle\dx\dt  + c\int_0^\tau \mathcal D(t)\dt\\
 -\int_{Q_\tau} \int_\Omega  \Big( \langle\nu_{t,x}, s|\vv|^2\rangle \rho(t,y) - \langle \nu_{t,x}, s\vv\rangle\langle \nu_{t,y}, s\vv \rangle\Big)\psi(x-y)\dy \dx \dt,
\end{multline*}
holds for  a.a. $\tau \in (0,T)$.  
\end{itemize}
With   ${\nu_0}_x = \delta_{(\rho_0(x) , \u_0(x))}$ for a.a. $x\in \Omega$ and $c$ independent of $\varepsilon$ and any solution.
\end{definition}
\begin{remark}\label{rem.post}
It is not difficult to see that 
\begin{multline*}
\int_\Omega \int_\Omega\Big( \langle \nu_{t,x}, s|\vv|^2\rangle \rho(t,y)  - \langle \nu_{t,x}, s\vv \rangle \langle \nu_{t,y}, s\vv \rangle\Big)\psi(x-y) \dy \dx \\ = \int_\Omega \int_\Omega \left(\frac 12 \langle \nu_{t,x}, s|\vv|^2 \rangle \rho(t,y)   -\langle \nu_{t,x}, s\vv \rangle\langle \nu_{t,y}, s\vv \rangle\right.\\ \left. +  \frac 12 \rho(t,x) \langle \nu_{t,y}, s|\vv|^2 \rangle\right) \psi(x-y)\dy \dx
=:\int_\Omega \int_\Omega \mathcal I \dy \dx,
\end{multline*}
is positive by the use of the Fubini lemma, 
\begin{multline*}
\mathcal I = \psi(x-y)\left(\frac 12\int s|\vv|^2\ {\rm d}\nu_{t,x}(s,\vv) \int \sigma\ {\rm d}\nu_{t,y}(\sigma, {\bf w})\right. \\ \left.- \int s\vv\ {\rm d}\nu_{t,x}(s,v) \int \sigma {\bf w} \ {\rm d}\nu_{t,y}(\sigma, {\bf w})
+ \frac 12 \int s\ {\rm d} \nu_{t,x}(s,\vv)\int \sigma {\bf w}^2\ {\rm d} \nu_{t,y}(\sigma, {\bf w})\right)\\
 = \frac 12 \int \int s\sigma |\vv-{\bf w}|^2 \ {\rm d}\nu_{t,x}(s,\vv)\ {\rm d} \nu_{t,y}(\sigma, {\bf w})\geq 0.
\end{multline*}
Here, the integration is over $[0,\infty)\times  \mathbb R^N$.

\end{remark}
In this paper we show that the concept of dissipative measure-valued solution is an appropriate one for our problem. On the one hand it is general enough to allow for existence, otherwise unavailable for weak solutions, and for another it is robust enough to provide an inviscid limit. The existence of measure-valued solutions is given as follows.

\begin{theorem}
\label{thm.exists}Let $\varepsilon >0$ and $\rho_0^\vep$, $\u_0^\vep$ be initial data with a finite energy, i.e.,
\begin{equation} \label{init.data.fin.en} \int_\Omega \left( \frac12 \rho_0^\vep |\u^\vep_0|^2 + P (\rho^\vep_0) \right)dx < \infty.\end{equation}
 If the hypotheses  \eqref{tlak.pred}, \eqref{h.pred} and \eqref{k.psi.pred} are satisfied, then for any $T> 0$ there exists a dissipative measure-valued solution  to \eqref{art.vis} and  \eqref{art.vis.init.con} in the sense of Definition \ref{dis.mes.sol}.
\end{theorem}

As already mentioned, we consider the dissipative measure-valued solutions to accommodate for the  nonlinearities in \eqref{art.vis}$_2$. We would like to emphasize that when $p(\rho)\sim \rho^\gamma,\ \gamma >\frac 32$ and $H(z)\equiv const.$ one can obtain the existence of weak solutions in a  standard manner (cf. \cite{Fe}). 

Before stating the main result  we would like to introduce the so-called {\it  relative entropy functional}. Let $\nu:\Omega \to \mathcal P ([0,\infty) \times \mathbb R^N)$ be a Young measure and let $r:\Omega\mapsto  (0,\infty)$ and $\U:\Omega \mapsto \mathbb R^N$ be smooth functions.
We then define the relative entropy functional $\mathcal E$ as
$$
\mathcal E(\nu, r, \U):= \int_\Omega \left\langle\nu_{x}, \frac 12 s |\vv-\U(x)|^2 + P(s) - P'(r(x))(s-r(x)) - P(r(x)) \right\rangle \dx.
$$

\begin{theorem}\label{thm.vis.lim}
Let $T> 0$ and  the hypotheses  \eqref{tlak.pred}, \eqref{h.pred} and \eqref{k.psi.pred} are satisfied. Let $\{(\nu^\vep,D^\vep)\}_{\vep>0}$ be a family of dissipative measure-valued solutions to \eqref{art.vis}  emanating from the initial conditions \eqref{art.vis.init.con}. Let $(r,\U)$ be a strong solution to \eqref{zaklad} in $Q_T$ endowed with \eqref{init.con} such that $r\geq c>0$ for some $c$. Moreover, assume that 
$$(\rho^\vep_0,\u^\vep_0)\rightarrow (r_0,\U_0),$$
in the following  sense:
$$
\mathcal E(\nu^\vep_0,r_0,\U_0)\rightarrow 0,
$$
where ${{\nu^\vep_0}_x} = \delta_{(\rho^\vep_0(x),\u_0^\vep(x))}$ for a.a $x\in \Omega$.

Then $\{(\nu^\vep,D^\vep)\}_{\vep>0}$ tends to $(r,\U)$ in the following sense:
\begin{equation}\label{main.res}
\underset{{\tau\in(0,T)}}\esssup \mathcal E(\nu^\vep_{\tau,\cdot},r(\tau),\U(\tau)) \rightarrow 0\ \ \mbox{ and } \ \  \underset{{\tau\in(0,T)}}\esssup \mathcal {D}^\vep(\tau) \rightarrow 0,
\end{equation}
as $\varepsilon \to 0$.
\end{theorem}
\begin{remark} Due to the assumptions there exist $\underline r,\overline r>0$ such that $\underline r\leq \frac {r(t,x)}2 \leq 2r(t,x) \leq \overline r$ for all $(t,x)\in Q_T$. According to \cite[(4.1)]{BuFeNo} we have
\begin{multline*}
\int_{(0,\infty)\times \mathbb R^N} |s-r(\tau,x)|\ {\rm d}\nu^\vep_{\tau,x}(s,\vv) \leq  \left(\int_{[\underline r,\overline r]\times \mathbb R^N} |s-r(\tau,x)|^2\ {\rm d}\nu^\vep_{\tau,x}(s,\vv) \right)^{\frac 12}\\ + c \left[\left(\int_{ (\overline r,\infty)\times \mathbb R^N}s \ {\rm d}\nu^\vep_{\tau,x}(s,\vv) \right)^{\frac 12}  + \left(\int_{(0,\underline r)\times \mathbb R^N} 1 \ {\rm d}\nu^\vep_{\tau,x} (s,\vv)\right)^{\frac 12}\right]^2\\  \leq c \langle \nu^\vep_{\tau,x} , P(s) - P'(r(\tau,x))(s-r(\tau,x)) - P(r(\tau,x))\rangle^{\frac 12} \left( 1 + \rho^\vep(\tau,x)\right)^{\frac 12},
\end{multline*}
and thus by H\" older inequality we get
\begin{multline}\label{odhad.e}
\|(\rho^\vep - r)(\tau)\|_{L^1(\Omega)} \leq \int_\Omega \langle \nu^\vep_{\tau,x} ,|s-r(\tau,x)|\rangle\dx \\\leq \int_{\Omega} \langle \nu^\vep_{\tau,x} , P(s) - P'(r(\tau,x))(s-r(\tau,x))-P(r(\tau,x))\rangle^{\frac 12}(1+\rho^\vep(\tau,x))^{\frac 12}\dx\\ \leq c\mathcal E(\nu^\vep_{\tau,\cdot} ,r(\tau),\U(\tau))^{\frac 12},
\end{multline}
where we use the fact that $\int_\Omega \rho^\vep(\tau)\dx =\int_\Omega \rho^\vep_0\dx$ for all $\tau\in [0,T]$ and $\|\rho_0^\vep \|_{L^1(\Omega)} \leq c$ for some $c>0$ independent of $\vep$ and any solution.

Thus \eqref{main.res} implies that
$$
\rho^\vep \rightarrow r\quad \mbox{in }L^\infty(0,T,L^1(\Omega)),
$$
as $\varepsilon \to 0$.
\end{remark}

It is worth pointing out that  the result of Theorem \ref{thm.vis.lim} applies also to suitable weak solutions (see \cite{BuFeNo}) as every suitable weak solution can be seen as a dissipative measure-valued solution.

The rest of this paper is devoted to the proofs of Theorems \ref{thm.exists} and \ref{thm.vis.lim}. In Section \ref{exist.mv} we show the existence of dissipative measure-valued solutions. Section \ref{relenin} is concerned with {\it relative entropy inequality} for the dissipative measure-valued solutions. Finally, this is used in Section \ref{invis.limit} to prove Theorem \ref{thm.vis.lim}.

\section{Existence of measure-valued solutions}\label{exist.mv}
Throughout this section let $T>0$ and $\varepsilon >0$ be fixed and we omit the index $\vep$ for the sake of clarity. We construct a measure-valued solution to \eqref{art.vis} and \eqref{art.vis.init.con} by means of a two-level approximation. First, we introduce a system with viscous penalization of the density. Second, we solve it by the Faedo-Galerkin approximation. However, as we would like to proceed  in both limits at the same time, we take  the coefficient of artificial viscosity dependent on $n$ -- the coefficient of the Faedo-Galerkin approximation. 

\subsection{Faedo-Galerkin approximation}
We consider a family of nested  finite dimensional spaces $X_n$, $n=1,2,\dots$ consisting of smooth vector-valued functions defined on $\Omega$.

The equation of continuity  \eqref{art.vis}$_1$ is regularized using vanishing artificial  viscosity  as 
\begin{equation}\label{van.vis.reg}
\begin{split}\partial_t \rho^n + \diver_x(\rho^n \u^n) &= \frac 1n \Delta_x \rho^n \ \ \ \mbox{ in } Q_T,\\
\rho^n(0) &= \rho^n_{0} \  \ \ \ \ \ \ \ \ \,  \mbox{ in } \Omega.
\end{split}
\end{equation}

We continue in the spirit of \cite[Chapter 7]{Fe} and we look for the Faedo-Galerkin approximate solutions, namely, we look for velocities $\u^n \in C([0,T];X_n)$ that satisfy the integral identity
\begin{multline}\label{FGeqmomentum} \int_\Omega \rho^n \u^n(\tau)\cdot \boldsymbol\eta \dx  -  \int_\Omega {\bf m}^n_{0}\cdot \boldsymbol\eta \dx\\
=\int_{Q_\tau}  (\rho^n \u^n\otimes \u^n -\varepsilon \nabla_x \u^n):\nabla_x \boldsymbol\eta + p(\rho^n)\diver_x \boldsymbol \eta -\frac 1n (\u^n \otimes\nabla_x \rho^n):\nabla_x \boldsymbol \eta \dx \dt\\
+  \int_{Q_\tau} \left((1-H(|\u^n|^2))\rho^n \u^n\cdot \boldsymbol\eta - \rho^n \int_{\Omega} \nabla_xK(x-y)\rho^n(y) \dy\cdot \boldsymbol\eta \right. \\ \left.+ \rho^n \int_{\Omega} \psi(x-y) \rho^n(y) \left(\u^n(y)-\u^n(x)\right) \dy\cdot \boldsymbol\eta\right) \dx\dt,\end{multline}
for any test function $\boldsymbol \eta \in X_n$ and all $\tau \in [0,T]$, where $\rho^n = S_{\rho_0^n}(\u^n)$ is a unique solution to \eqref{van.vis.reg}.

As the additional terms on the right-hand side of \eqref{FGeqmomentum} can be treated in a standard way, we may  directly use the method\footnote{The only difference is the energy inequality which can be replaced by \eqref{ei.for.n}.} of \cite[Chapter 7]{Fe} to claim the existence of a unique solution to  \eqref{FGeqmomentum}  for any initial data that satisfy:
$$\rho^n_0 \in C^{2+ \delta}(\Omega),  \delta >0, \ \ \ \ \inf_\Omega \rho_0^n >0\ \  \ \mbox{ and }\ \ \ {\bf m}_0^n  \in L^2(\Omega).$$ 

\subsection{Limit in $n$}
The aim of this section is to tend with $n$ to infinity.
In order to derive the energy estimates we use $\boldsymbol \eta = \u^n$ as a test function in \eqref{FGeqmomentum} which together with \eqref{van.vis.reg} and the symmetry of $K$ and $\psi$ yield 
\begin{multline}\label{ei.for.n}
 \int_\Omega\left( \frac 12 \rho^n |\u^n|^2 + P(\rho^n) + \frac 12\rho^nK*\rho^n \right)(\tau)\dx\\ + \int_{Q_\tau} \left( \varepsilon |\nabla_x \u^n|^2 + \frac 1n P''(\rho^n)|\nabla_x \rho^n|^2\right) \dx\dt  \\ \leq \int_\Omega\left( \frac 12 \frac{|{\bf m}^n_0|^2}{\rho_0^n} + P(\rho_0^n) + \frac 12\rho_0^n K*\rho_0^n \right)\dx\\ + \int_{Q_\tau} (1-H(|\u^n|^2)) \rho^n |\u^n|^2\dx\dt + \frac 1n \int_{Q_\tau} \int_\Omega \rho^n(x)\Delta_x K(x-y)\rho^n(y) \dy \dx\dt \\
 -\frac 12 \int_{Q_\tau} \int_\Omega  \rho^n(x)\rho^n(y) \psi(x-y) \left(\u^n(x) - \u^n(y)\right)^2 \dy \dx \dt,
\end{multline}
for all $\tau \in [0,T]$.
From now on,  we shall consider $\rho_0^n$ and ${\bf m}_0^n$ such that $\rho_0^n \to \rho_0$ and ${\bf m}_0^n \to \rho_0 \u_0$ in $L^1(\Omega)$ together with 
$$\int_\Omega \left(\frac12 \frac{ |{\bf m}_0^n|^2}{\rho_0^n} + P(\rho_0^n) \right)\dx \to \int_\Omega \left(\frac12 \rho_0 |\u_0|^2 + P(\rho_0)\right)\dx.$$ This possible due to  \eqref{init.data.fin.en} and \cite[Section 7.10.7]{NoSt}.
 
By integrating \eqref{van.vis.reg} over $\Omega$ we deduce that 
\begin{equation} \label{est.rho.n}
\underset{t\in (0,T)}{\sup} \|\rho^n(t)\|_{L^1(\Omega)} <c,
\end{equation}
and since 
$$
 \frac1n\left|\int_\Omega \int_\Omega \rho^n(x)\Delta_x K(x-y)\rho^n(y) \dy \dx\right|\leq \frac1n \|\Delta_x K\|_{L^\infty(\Omega)} \|\rho^n(t)\|_{L^1(\Omega)}^2\leq \frac{c}n\to 0,
$$
we infer from \eqref{ei.for.n} and Gronwall inequality that
\begin{equation}\label{apriori.est}
\begin{split}
\underset{t\in (0,T)}\sup \|\rho^n(t)|\u^n|^2(t)\|_{L^1(\Omega)} &\leq c,\\
\underset{t\in (0,T)}\sup \|P(\rho^n(t))\|_{L^1(\Omega)} & \leq c,\\
\vep\|\nabla_x \u^n\|^2_{L^2(Q_T)} &\leq c,\\
\frac 1n \left\|P''(\rho^n)|\nabla_x \rho^n|^2\right\|_{L^1(Q_T)} & \leq c,
\end{split}
\end{equation}
with $c>0$ independent of $n$. Here and hereafter in this section, $c$ denotes a positive generic constant that is independent of $\vep$ and $n$ unless specified otherwise.  

From \eqref{tlak.pred} and \eqref{apriori.est}$_4$ we derive
\begin{equation}\label{apriori.est2}
\begin{split}
\underset{t\in (0,T)}\sup\int_\Omega \rho^n \log \rho^n \dx &\leq c,\\
\frac 1n \int_0^T\int_{\{\rho^n > 1\}} \frac{|\nabla_x\rho^n|^2}{\rho^n} \dx \dt &\leq c.
\end{split}
\end{equation}
As $\rho^n$ is smooth enough, we may multiply $\eqref{van.vis.reg}$ by $b'(\rho^n)$ for any $b\in C_c^{\infty}([0,\infty))$ in order to get 
\begin{multline*}
\pat b(\rho^n) + \diver_x (b(\rho^n)\u^n) + (b'(\rho^n)\rho^n - b(\rho^n)) \diver_x \u^n \\ = \frac 1n \diver_x (b'(\rho^n)\nabla_x\rho^n)  - \frac 1n b''(\rho^n)|\nabla_x \rho^n|^2 \mbox{ in } Q_T.
\end{multline*}
By choosing a suitable $b$, namely $b(z) = z^2$ for $z\leq 1$, and integrating over $Q_T$ we get
\begin{equation}\label{apriori.est3}
\frac 1n \int_0^T\int_{\{\rho^n\leq 1\}}|\nabla_x \rho^n|^2 \dx \dt \leq c(\vep).
\end{equation}

Next we obtain relative compactness in $L^1$. First, we deduce from \eqref{apriori.est2}$_1$  that $\{\rho^n\}$ is equi-integrable. Indeed, we have
$$
\log k\int_{\{\rho^n\geq k\}}\rho^n \dx\leq \int_{\{\rho^n\geq k\}}\rho^n \log \rho^n  \dx\leq c,\quad \forall k\in (0,\infty),
$$
yielding that $\underset{n}{\sup}\int_{\{\rho^n\geq k\}}\rho^n \dx\to 0$ as $k\to \infty$. 
Second, the sequence $\{\rho^n \u^n\}$ is also equi-integrable. To see this it is enough to consider the following sequence of inequalities, which holds for every $k,m\in (0,\infty)$,
\begin{multline*}
\int_{\{|\rho^n \u^n|\geq k\}}|\rho^n \u^n|\dx\leq \int_{\{\rho^n \geq \frac km\}} |\rho^n \u^n| \dx + \int_{\{|\u^n|\geq m\}} |\rho^n \u^n |\dx\\ 
\leq \frac{1}{\sqrt{\log \frac km}} \int_\Omega \sqrt{\rho^n\log \rho^n}\sqrt \rho^n \u^n \dx+ \frac 1m \int_\Omega \rho^n |\u^n|^2\dx\leq c\left(\frac1{\sqrt{\log \frac km}} + \frac 1m\right).
\end{multline*}
As the previous estimate holds for every choice of $k$ and $m$, it suffices to consider $m = \sqrt k$ in order to deduce 
$$
\sup_n \int_{\{|\rho^n\u^n|>k\}} |\rho^n \u^n| \dx \to 0,
$$
 as $k\to \infty$. Third, we can claim the equi-integrability of $\{\u^n\}$ using the Poincar\' e inequality which can be found, e.g., in \cite[Lemma 3.2]{Fe}.

Now we are able to proceed to limit in a weak formulation of the continuity equation \eqref{van.vis.reg}, i.e., in
\begin{multline} \label{van.vis.weak}
\int_\Omega \rho^n(\tau) \varphi(\tau)\dx - \int_\Omega \rho_0^n\varphi (0) \dx - \int_{Q_\tau} \rho^n \pat \varphi \dx \dt  - \int_{Q_\tau}\rho^n \u^n  \nabla_x \varphi\dx \dt \\= - \frac1n \int_{Q_\tau} \nabla_x \rho^n \nabla_x \varphi \dx \dt ,
\end{multline}
which holds for every $\varphi \in C^1([0,T],\Omega)$ and $\tau \in [0,T]$. By \eqref{est.rho.n}, \eqref{apriori.est2}$_2$ and \eqref{apriori.est3} the right-hand side of \eqref{van.vis.weak} vanishes  as $n\to \infty$ since
\begin{multline*}
\left|\frac1n\int_0^T\int_{ \{\rho^n> 1\}} \frac 1{\sqrt{\rho^n}} \nabla_x \rho^n\sqrt{\rho^n} \nabla_x \varphi\dx \dt + \frac 1n \int_0^T\int_{\{\rho^n \leq 1\}} \nabla_x \rho^n\nabla_x \varphi \dx \dt\right|\\ \leq \frac 1{\sqrt n} \left(\frac1n\int_0^T\int_{\{\rho^n> 1\}} \frac{|\nabla_x \rho^n|^2}{\rho^n} \dx \dt\right)^{\frac 12} \|\nabla_x \varphi\|_{L^\infty(\Omega)} \|\rho^n\|_{L^1(Q_T)}^{\frac 12} \\ + \frac 1{\sqrt n} \left( \frac1n\int_0^T\int_{\{\rho^n \leq 1\}} |\nabla_x \rho^n|^2 \dx \dt\right)^\frac12 \|\nabla_x \varphi\|_{L^2(Q_T)} \leq \frac 1{\sqrt n} c(\vep).
\end{multline*}
Next we may use \cite[Theorem 6.2]{Ped} and there exists a Young measure 
$$\nu_{t,x} \in \mathcal P([0,\infty)\times \mathbb R^N) \mbox{ for a.a. }(t,x) \in [0,T]\times \Omega,$$ 
associated with the equi-integrable sequence $\{\rho^n, \u^n\}$. We will use the notation $\rho(t,x) = \langle \nu_{t,x},s\rangle$ and $\u(t,x)= \langle \nu_{t,x},\vv\rangle$ for shortness. As we proceed to limit in \eqref{van.vis.weak} we get
\begin{multline*}
\int_\Omega \rho(\tau) \varphi(\tau)\dx - \int_{Q_\tau} \rho  \pat \varphi \dx \dt  - \int_{Q_\tau}\langle \nu_{t,x}, s\vv\rangle \nabla_x \varphi(t,x)\dx \dt = \int_\Omega \rho_0\varphi (0) \dx,
\end{multline*}
which holds for every $\varphi \in C^1([0,T],\Omega)$ and all $\tau \in [0,T]$ as $\rho^n$ is precompact in $C_{weak} ([0,T], L^1(\Omega))$. 

Next, we show that we can proceed to limit in the balance of energy \eqref{ei.for.n}. Due to $L^1$-bounds from \eqref{apriori.est} we get
\begin{equation*}
\begin{split}
\frac 12 \rho^n |\u^n|^2 + P(\rho^n) &\to \left\langle \nu_{t,x}, \frac 12 s|\vv|^2 + P(s)\right\rangle + \xi\ \mbox{weakly}^* \mbox{ in }L^\infty_{\weak}(0,T,\mathcal M(\Omega)),\\
 |\nabla_x \u^n|^2 &\to |\nabla_x \u|^2 + \sigma \ \mbox{weakly}^* \mbox{ in } \mathcal M^+([0,T]\times\Omega),
\end{split}
\end{equation*}
for some non-negative measures $\xi$ and $\sigma$. Further, since $K$ is smooth and $\{\rho^n\}$ is equi-integrable we obtain 
$$
K*\rho^n \to  K*\rho\ \mbox{a.e. in }Q_T  \mbox{ and } \|K*\rho^n \|_{L^\infty(Q_T)}\leq c.
$$
We use Young inequality with the functions $\Phi(t) = (t+1)\log(t+1) -t$ and $\Psi(t) = e^t - t - 1$ in order to derive
$$
\int_{\Omega} \rho^n \left(K*\rho^n- K*\rho\right) \dx \leq \delta \int_{\Omega} \Phi(\rho^n) \dx + c(\delta)\int_\Omega \Psi\left(K*\rho^n - K*\rho \right)\dx,
$$
for arbitrary $\delta>0$. Since the second term on the right-hand side tends to zero as $n\to \infty$ due to  Lebesgue's dominated convergence theorem and the first term can be estimated by $\delta c$, we immediately conclude that
$$
\frac12\int_\Omega \left(\rho^n K*\rho^n\right)(\tau) \dx \rightarrow \frac12\int_\Omega (\rho K*\rho)(\tau) \dx,
$$
for all $\tau \in [0,T]$. Since $P''(z) = \frac{p'(z)}z > 0$ we neglect the term $\frac1n \int_{Q_\tau} P''(\rho^n) |\nabla_x \rho^n|^2\dx\dt$.
Next, we get
$$
\int_{Q_\tau} (1-H(|\u^n|^2))\rho^n|\u^n|^2\dx \dt  \rightarrow \int_{Q_\tau} \left\langle \nu_{t,x}, (1-H(|\vv|^2))s|\vv|^2\right\rangle \dx \dt  + \int_0^\tau \kappa(t)|\Omega|\dt,
$$
for all $\tau \in [0,T]$ and  some  measure $\kappa$. Due to \cite[Lemma 2.1]{FeGwSwWi} and boundedness of $H$  we have $|\kappa(t)|\leq c \xi(t)$ for all $t\in [0,T]$.

In order to derive the convergence in the last term of \eqref{ei.for.n}  it is enough to realize that the Arzela-Ascoli theorem yields for all $t\in [0,T]$ that
\begin{equation*}
\begin{split}
\psi*\rho^n &\rightrightarrows \psi*\rho \mbox{ on } \Omega,\\
\psi*(\rho^n\u^n) &\rightrightarrows \psi*\langle\nu_{t,y},s\vv\rangle \mbox{ on } \Omega,
\end{split}
\end{equation*}
as far as $\psi\in C^1(\Omega)$. Thus we may infer that 
\begin{multline*}
\lim_{n\to \infty}\left(\int_\Omega \int_\Omega \rho^n(x)\rho^n(y)\psi(x-y) |\u^n(x)|^2\dy\dx\right.\\ \left.- \int_\Omega \int_\Omega\rho^n(x)\rho^n(y) \psi(x-y) \u^n(x)\u^n(y))\dy \dx\right)(t)\\
= \int_\Omega \int_\Omega \langle\nu_{t,x}, s|\vv|^2\rangle \rho(t,y) \psi(x-y)\dy \dx \\- \int_\Omega \int_\Omega \langle \nu_{t,x}, s\vv\rangle\psi(x-y)\langle \nu_{t,y}, s\vv \rangle \dy \dx + \int_\Omega \psi*\rho(t) {\rm d}\lambda(t),
\end{multline*}
where $\lambda$ is a non-negative measure satisfying $\lambda(t) \leq c \xi(t)$ for all $t\in [0,T].$

As a consequence of the previous ideas we get that \eqref{ei.for.n} converges to 
\begin{multline*}
\int_\Omega \left\langle\nu_{\tau,x}, \frac 12s|\vv|^2 + P(s) \right\rangle\dx + \frac 12\int_\Omega (\rho K*\rho)(\tau)\dx \\
+ \vep  \int_{Q_\tau}  |\nabla_x \u|^2\dx \dt + \xi(\tau)|\Omega|  \\
 \leq\int_{\Omega} \left( \frac 12 \rho_0 |\u_0|^2 + P(\rho_0) + \frac 12 \rho_0(K*\rho_0)\right) \dx \\
+\int_{Q_\tau} \langle\nu_{t,x}, (1- H(|\vv|^2)) s|\vv|^2\rangle\dx\dt  + c \int_0^\tau \xi(t)|\Omega| \dt\\
 -\int_{Q_\tau}  \int_\Omega \Big( \langle\nu_{t, x}, s|\vv|^2\rangle \rho(t,y) \psi(x-y)- \langle \nu_{t,x}, s\vv\rangle\psi(x-y)\langle \nu_{t,y}, s\vv \rangle\Big)\dy \dx\dt,
\end{multline*}
for a.a. $\tau \in (0,T)$.

It remains to proceed to limit in the momentum equation \eqref{FGeqmomentum}. Nevertheless, the term $\frac 1n \int_{Q_\tau} (\u^n \otimes \nabla_x \rho^n): \nabla_x \bphi \dx \dt$ vanishes due to $L^2$ estimates following from \eqref{apriori.est}$_1$, \eqref{apriori.est2}$_2$, \eqref{apriori.est3} and Poincar\' e inequality and all the other terms are equi-integrable except of $\int_{Q_\tau} p(\rho^n)\diver_x \bphi \dx \dt$ and $\int_{Q_\tau} (\rho^n \u^n \otimes \u^n ):\nabla_x \bphi \dx \dt$. However, since $|\rho^n \u^n\otimes \u^n|\leq \rho^n|\u^n|^2$ and $p(\rho^n)\leq a P(\rho^n)$ for $\rho^n$ sufficiently large and  some $a>0$, we get
\begin{multline*}
\int_\Omega \langle \nu_{\tau,x}, s\vv\rangle \bphi(\tau,x) \dx  - \int_{Q_\tau}\langle \nu_{t,x},s\vv\rangle \pat\bphi(t,x)\dx \dt\\
 - \int_{Q_\tau} \langle \nu_{t,x}, s\vv\otimes \vv\rangle \nabla_x \bphi(t,x)\dx \dt + \vep\int_{Q_\tau}\nabla_x \u : \nabla_x\bphi \dx \dt\\ - \int_{Q_\tau} \langle \nu_{t,x},p(s)\rangle \diver_x \bphi(x) \dx \dt  =\int_{Q_\tau} \langle \nu_{t,x},  (1-H(|\vv|^2))s\vv\rangle \bphi(t,x)\dx \dt \\- \int_{Q_\tau} \rho(t,x) \int_\Omega \nabla_xK(x-y) \rho(t,y) \dy \bphi(t,x)\dx \dt\\ 
+ \int_{Q_\tau} \rho(t,x)\int_\Omega \psi(x-y)\langle \nu_{t,y},s\vv \rangle\dy \bphi(t,x)\dx \dt\\
 - \int_{Q_\tau} \langle \nu_{t,x}, s\vv\rangle \int_\Omega \psi(x-y)\rho(t,y)\dy \bphi(t,x)\dx \dt\\
 + \int_\Omega \rho_0 \u_0 \bphi(0) \dx + \int_0^\tau \int_\Omega \nabla_x \bphi(t) {\rm d}\mu(t) \dt,
\end{multline*}
where, due to \cite[Lemma 2.1]{FeGwSwWi}, $\mu$ is a measure such that $|\mu(t)|\leq c \xi(t)$ for all $t\in[0,T]$.

This concludes the proof of Theorem \ref{thm.exists}.

\section{Relative entropy inequality}\label{relenin}
Let $(\nu^\varepsilon, \mathcal D^\vep)$ be a measure-valued solution to \eqref{art.vis} and  \eqref{art.vis.init.con}. Similarly as in \cite[Section 3]{FeGwSwWi} one can show that every pair of smooth functions $(r,\U)\in C^1([0,T]\times \Omega)\times C^1([0,T]\times \Omega)^N$ satisfies
\begin{multline}
\label{REI}\mathcal E(\nu^\vep_{\tau,\cdot}, r(\tau), \U(\tau)) + \mathcal D_\vep(\tau) + \int_\Omega\frac 12\Big((\rho^\vep - r) K*(\rho^\vep- r)\Big)(\tau)\dx  \\ + \varepsilon  \int_{Q_\tau}\nabla_x (\u^\vep - \U):\nabla_x(\u^\vep - \U) \dx \dt 
 \leq \mathcal E(\nu^\vep_0, r(0),\U(0))\\ + \int_\Omega\frac 12\Big((\rho^\vep - r) K*(\rho^\vep- r)\Big)(0)\dx +  \mathcal R(\nu^\vep_{\tau,\cdot}, r(\tau), \U(\tau)), 
\end{multline}
where
\begin{multline}\label{REIrest}
\mathcal R(\nu^\vep_{\tau,\cdot}, r(\tau), \U(\tau)) = \int_{Q_\tau} \Big( \langle \nu^\vep_{t,x}, s(\U-\vv)\rangle \pat \U + \langle \nu^\vep_{t,x}, s\vv\otimes(\U - \vv)\rangle\nabla_x \U\Big)\ {\rm d}x \ {\rm d}t\\  +\vep\int_{Q_\tau} \nabla_x \U :\nabla_x(\U-\u^\vep) \ {\rm d}x \ {\rm d}t \\+ 
\int_{Q_\tau}\left(\left(1-\frac{\rho^\vep}r\right) p'(r)\pat r - \langle \nu^\vep_{t,x}, s\vv\rangle \frac{p'(r)}{r}\nabla_x r + \langle \nu^\vep_{t,x}, p(s)\rangle \diver_x \U\right)\ {\rm d}x \ {\rm d}t \\ 
+  \int_{Q_\tau} \langle \nu_{t,x}^\vep, s\vv(\vv - \U)(1-H(|\vv|^2))\rangle \ {\rm d}x \ {\rm d}t 
+\int_{Q_\tau} (\rho^\vep \nabla_x K*\rho^\vep) \cdot \U \ {\rm d}x \ {\rm d}t\\ - \frac 12\int_{\Omega} \Big(\big(r K*\rho^\vep + (\rho^\vep - r)K*r\big)(\tau) - \big(r K*\rho^\vep + (\rho^\vep - r)K*r\big)(0)\Big) \dx\\
+  \int_{Q_\tau} \int_\Omega\langle \nu^\vep_{t,x}, s(\vv-\U(t,x))\rangle \langle\nu^\vep_{t,y}, s\vv\rangle\psi(x-y)\dy \dx \dt\\ - \int_{Q_\tau} \int_\Omega \langle \nu^\vep_{ t,x}, s\vv(\vv-\U(t,x))\rangle \rho^\vep(t,y) \psi(x-y)\dy \dx \dt\\
 + \int_0^\tau \Big( \int_\Omega \nabla_x \U {\rm d} \mu^\vep(t) + c \mathcal D^\varepsilon(t)\Big)\dt,
\end{multline}
for a.a. $\tau \in (0,T)$ and $c>0$ independent of $\varepsilon$ and any solution.

\section{Inviscid limit}\label{invis.limit}
Let $(r, \U)$ be a strong solution to \eqref{zaklad} and  \eqref{init.con}. Consequently,
\begin{multline*}
\int_{Q_\tau} (\rho^\vep \nabla_x K*\rho^\vep) \cdot \U \ {\rm d}x \ {\rm d}t \\= \int_{Q_\tau} \Big( (\rho^\vep-r )\nabla_x K* (\rho^\vep - r)\cdot \U + r\nabla_x K*(\rho^\vep - r) \cdot\U + (\rho^\vep \nabla_x K* r )\cdot\U \Big)\dx \dt\\
 = \int_{Q_\tau}  \Big( (\rho^\vep-r ) \U \cdot \nabla_x K* (\rho^\vep - r) + \pat r K*(\rho^\vep - r)  \\ + \langle \nu^\vep_{t,x}, s(\U-\vv)\rangle\nabla_x K*r + \pat \rho^\vep K*r \Big)\ {\rm d}x \ {\rm d}t.
\end{multline*}
%
We rewrite \eqref{REI} and \eqref{REIrest} in the following way
\begin{multline}
\label{REIsUr}
\mathcal E(\nu^\vep_{\tau,\cdot}, r(\tau), \U(\tau)) + \vep \int_{Q_\tau} |\nabla_x (\u^\vep - \U)|^2  \dx \dt + \mathcal D^\varepsilon (\tau)\\
 +\frac12 \int_\Omega\Big( \big((\rho^\vep - r) K*(\rho^\vep- r)\big)(\tau) -\big((\rho^\vep - r) K*(\rho^\vep- r)\big)(0)\Big)\dx\\
\leq \mathcal E(\nu^\vep_0, r_0,\U_0)\\
 +\int_{Q_\tau} \Big( \langle \nu^\vep_{t,x}, s(\U-\vv)\rangle \pat \U + \langle \nu^\vep_{t,x}, s\vv\otimes(\U - \vv)\rangle\nabla_x \U \Big)\ {\rm d}x \ {\rm d}t\\
  +\vep\int_{Q_\tau} \nabla_x \U :\nabla_x(\U-\u^\vep) \ {\rm d}x \ {\rm d}t \\
+\int_{Q_\tau} \langle \nu^\vep_{t,x}, s\U - s\vv\rangle \frac{p'(r)}r\nabla_x r \dx \dt\\
-\int_{Q_\tau} \langle \nu^\vep_{t,x},  p(s) - p'(r)(s-r) - p(r)\rangle\diver_x \U \dx \dt\\
+  \int_{Q_\tau} \langle \nu^\vep_{t,x}, s\vv(\vv - \U)(1-H(|\vv|^2))\rangle \ {\rm d}x \ {\rm d}t\\
+\int_{Q_\tau} \Big((\rho^\vep - r) \U \nabla_x K*(\rho^\vep - r) + \langle \nu^\vep_{ t, x}, s(\U - \vv) \rangle \nabla_x K* r \Big)\dx \dt\\
+  \int_{Q_\tau} \int_\Omega\langle \nu^\vep_{t,x}, s(\vv-\U(t,x))\rangle \langle\nu^\vep_{ t,y}, s\vv\rangle\psi(x-y)\dy \dx \dt\\ - \int_{Q_\tau} \int_\Omega \langle \nu^\vep_{ t,x}, s\vv(\vv-\U(t,x))\rangle \rho^\vep(t,y) \psi(x-y)\dy \dx \dt\\
+ \int_0^\tau \Big(\int_\Omega \nabla_x \U{\rm d}\mu^\vep(t) + c \mathcal D^\varepsilon(t)\Big)\dt =:\sum_{i=1}^{10}\mathcal I_i,
\end{multline}
for a.a. $\tau \in (0,T)$.
Here and hereafter $c$ denotes a positive generic constant independent of $\varepsilon$.
\subsection{Estimates of $\mathcal I_i$}
Since we assume that $\mathcal E (\nu^\vep_{0}, r_0, \U_0) \to 0$ as $\vep \to 0$, we get
\begin{equation*}
\mathcal I_1 \leq \Gamma(\vep),
\end{equation*}
where $\Gamma(\vep)\to 0$ as $\vep \to 0$.

Further, as $\U$ is smooth, we deduce
\begin{equation*}
\mathcal I_3 \leq \vep\frac { c}2+ \vep\frac12 \int_{Q_\tau} |\nabla_x (\u_{\vep} - \U)|^2\dx\dt.
\end{equation*}

According to estimates on $\mu^\vep $ and smoothness of $\U$, we get
\begin{equation*}
\mathcal I_{10} \leq c \int_0^\tau\mathcal D^\vep (t)\dt.
\end{equation*}

Due to \eqref{zaklad}$_1$
\begin{multline*}
\mathcal I_2 + \mathcal I_4 + \mathcal I_7 = \int_{Q_\tau} \bigg[ \langle \nu^\vep_{t,x}, s(\U - \vv)\rangle \cdot\\
\cdot\left(-\U \nabla_x \U + (1-H(|\U|^2)) \U + \int_\Omega\psi(x-y)r(y)(\U(y) - \U(x))\dy\right)\\
+ \langle \nu^\vep_{t,x}, s\vv\otimes (\U-\vv)\rangle \nabla_x \U + (\rho^\vep - r)\U \nabla_x K*(\rho^\vep -r)\bigg]\dx \dt\\
 = \int_{Q_\tau} \langle \nu^\vep_{t,x}, s(\vv-\U)\nabla_x \U (\U-\vv) \rangle \dx \dt  \\ 
 + \int_{Q_\tau} \langle \nu^\vep_{t,x}, s(\U - \vv)\rangle (1-H(|\U|^2)) \U\dx \dt\\ + \int_{Q_\tau} \langle \nu^\vep_{t,x}, s(\U - \vv)\rangle\int_\Omega\psi(x-y)r(y)(\U(y) - \U(x))\dy\dx \dt\\
 + \int_{Q_\tau}(\rho^\vep - r)\U \nabla_x K*(\rho^\vep -r)\dx \dt\\
=:\mathcal I_{11} +\mathcal I_{12} + \mathcal I_{13} + \mathcal I_{14}.
\end{multline*}
The Taylor formula and \cite[(4.1)]{BuFeNo} yield
$$
|p(s) - p'(r)(s-r) - p(r)|\leq c |P(s) - P'(r)(s-r) - P(r)|,
$$
which together with the definition of $\mathcal E$ imply
\begin{equation*}
\mathcal I_5 + \mathcal I_{11}\leq c\int_0^\tau\mathcal E (\nu^\vep_{t,\cdot}, r(t), \U(t))\dt.
\end{equation*}
Further,
$$
\mathcal I_6 + \mathcal I_{12} = \int_{Q_\tau} \langle \nu^\vep_{t, x}, s(\vv-\U) \Big(\vv(1-H(|\vv|^2)) - \U(1-H(|\U|^2))\Big)\rangle\dx \dt.
$$
Due to the assumptions   on $H$ it is possible to split it into two parts, namely $H = H_1 + H_2$, where $H_1\in C^1([0,\infty))$ has a compact support and $H_2$ is non-decreasing. It then  follows that 
$$
 - s(\vv-\U) \left(vH_2(|\vv|^2) - \U H_2(|\U|^2)\right)\leq 0.
$$
Since $H_1$ is a Lipschitz function with  compact support, we get $$|H_1(|\vv|^2) -H_1(|\U|^2)|\leq c |\vv-\U|,$$ and thus
\begin{multline*}
\langle \nu^\vep_{t,x}, s(\vv-\U)(\vv H_1(|\vv|^2) - \U H_1(|\U|^2))\rangle\\ \leq \langle \nu^\vep_{t,x} s(\vv-\U)^2 H_1(|\vv|^2)\rangle + c\langle \nu^\vep_{t,x}, s (\vv-\U)^2\rangle.
\end{multline*}
Consequently,
\begin{equation*}
\mathcal I_6 + \mathcal I_{12}\leq \int_0^\tau \mathcal E (\nu^\vep_{t,\cdot}, r(t), \U(t)) \dt.
\end{equation*}
Further, 
\begin{multline*}
\mathcal I_8 + \mathcal I_9 + \mathcal I_{13} 
= \int_{Q_\tau}\int_\Omega \psi(x-y) \langle \nu^\vep_{t,x}, s(\vv-\U)\rangle \langle \nu^\vep_{t,y}, (s-r)(\U - \U(x))\rangle \dy \dx\dt\\
- \frac 12 \int_{Q_\tau}\int_\Omega \psi(x-y)\Big(\langle\nu^\vep_{t,x}, s(\vv-\U)^2\rangle \rho^\vep(y) - \langle \nu^\vep_{t,x}, s(\vv-\U)\rangle \langle \nu_{t,y}^\vep, s(\vv-\U)\rangle\Big)\dy\dx\dt.
\end{multline*}
The second term on the right-hand side is negative by the same argument as in Remark \ref{rem.post}. By the Young's inequality we deduce that
$$
\left|\int_{(0,\infty)\times \mathbb R^N} s(\vv-\U) f(t,x,y)d\nu_{t,x}^\vep\right|\leq \frac 12 \int_{(0,\infty)\times \mathbb R^N} s(\vv-\U)^2 d\nu_{t,x}^\vep + \frac 12 \rho^\vep f^2(t,x,y),
$$
for any $f$ independent of $(s, \vv)$.
This together with \eqref{odhad.e} implies
\begin{multline*}
\mathcal I_8 + \mathcal I_9 + \mathcal I_{13}\leq \frac 12 \int_{Q_\tau} \langle \nu^\vep_{t,x}, s(\vv-\U)^2\rangle \dx \dt + c \int_0^\tau\|\rho^\vep - r\|_{L^1(\Omega)}^2 \dt \\ \leq c \int_0^\tau \mathcal E(\nu^\vep_{t,\cdot}, r(t),\U(t)) \dt.
\end{multline*}
Further, \eqref{odhad.e} immediately implies
\begin{equation*}
\mathcal I_{14} \leq c \int_0^\tau \mathcal E(\nu^\vep_{t,\cdot},r(t),\U(t)) dt.
\end{equation*}
Finally, by the same method as in \cite[Section 5]{CaFeGwSw} we derive
\begin{multline*}
\frac12 \int_\Omega\Big( \big((\rho^\vep - r) K*(\rho^\vep- r)\big)(\tau) -\big((\rho^\vep - r) K*(\rho^\vep- r)\big)(0)\Big)\dx\\
=  \int_{Q_\tau} \langle\nu^\vep_{ t,x}, r\U-s\vv \rangle \nabla_x K*(r-\rho^\vep )\dx\dt,
\end{multline*}
and
\begin{multline*}
\int_{Q_\tau} \langle\nu^\vep_{ t,x}, r\U-s\vv \rangle \nabla_x K*(r-\rho^\vep)\dx \dt = \\ \int_{Q_\tau} (r-\rho^\vep)\U \nabla_x K*(r-\rho^\vep)\dx \dt + \int_{Q_\tau} \langle \nu^\vep_{t,x}, s(\U - \vv)\rangle \nabla_x K*(r-\rho^\vep) \dx \dt\\ \leq c\int_0^\tau\mathcal E(\nu^\vep(t),r(t),\U(t))\dt + c\int_{Q_\tau} \langle \nu^\vep_{t,x},s|\U-\vv|^2\rangle \dx \dt \\+c \|\nabla_x K*(\rho^\vep - r)\|_\infty^2 \int_{Q_\tau} \rho^\vep \dx \dt \leq c\int_0^\tau\mathcal E(\nu^\vep_{t,\cdot},r(t),\U(t))dt,
\end{multline*}
where we have used the  H\"older inequality,
$$
\int_{(0,\infty)\times \mathbb R^N} \sqrt s \sqrt s(\U - \vv) \ {\rm d}\nu^\vep_{t,x}\leq \left(\int_{(0,\infty)\times \mathbb R^N} s \ {\rm d}\nu^\vep_{t,x}\right)^{\frac 12} \left(\int_{(0,\infty)\times \mathbb R^N} s|\U - \vv|^2 \ {\rm d}\nu^\vep_{t,x}\right)^{\frac 12}.
$$
\subsection{Limit}
We collect all the estimates from the previous subsection together with \eqref{REIsUr} in order to derive
$$
\mathcal E(\nu^\vep_{\tau,\cdot}, r(\tau), \U(\tau)) + \mathcal D^\vep(\tau) \leq c\int_0^\tau  \mathcal E(\nu^\vep_{t,\cdot},r(t), \U(t)) + \mathcal D^\vep(t) \dt + \Gamma(\vep),
$$
for a.a. $\tau \in (0,T)$ and, by the Gr\"onwall's inequality,
$$
\mathcal E(\nu^\vep_{\tau,\cdot}, r(\tau) , \U(\tau))  + \mathcal D^\vep(\tau)\leq c\Gamma(\vep),
$$ 
for a.a. $\tau\in (0,T)$, where $\Gamma(\vep) \to 0$ as $\vep \to 0$. This concludes the proof of Theorem \ref{thm.vis.lim}.

\section*{References}

\bibliography{bmbibfile}

\end{document}